\documentclass[11pt,reqno]{amsart}
\usepackage{amssymb,latexsym}
\usepackage{graphicx}
\usepackage{url}		

\calclayout

\makeatletter
\newcommand{\monthyear}[1]{%
  \def\@monthyear{\uppercase{#1}}}
\newcommand{\volnumber}[1]{%
  \def\@volnumber{\uppercase{#1}}}
\AtBeginDocument{%
\def\ps@plain{\ps@empty
  \def\@oddfoot{\@monthyear \hfil \thepage}%
  \def\@evenfoot{\thepage \hfil \@volnumber}}
\def\ps@firstpage{\ps@plain}
\def\ps@headings{\ps@empty
  \def\@evenhead{%
    \setTrue{runhead}%
    \def\thanks{\protect\thanks@warning}%
    \uppercase{The Fibonacci Quarterly}\hfil}%
  \def\@oddhead{%
    \setTrue{runhead}%
    \def\thanks{\protect\thanks@warning}%
    \hfill\uppercase{FIBONACCI NUMBERS WHICH ARE PRODUCTS OF
TWO PELL NUMBERS}}%
  \let\@mkboth\markboth
  \def\@evenfoot{%
    \thepage \hfil \@volnumber}%
  \def\@oddfoot{%
    \@monthyear \hfil \thepage}%
  }%
\footskip=25pt
\pagestyle{headings}%
}
\makeatother

\newcommand{\R}{{\mathbb R}}

\newcommand{\Z}{{\mathbb Z}}
\def \L {{\mathbb L}}

\theoremstyle{plain}
\numberwithin{equation}{section}
\newtheorem{thm}{Theorem}[section]
\newtheorem{theorem}[thm]{Theorem}
\newtheorem{lemma}[thm]{Lemma}

\begin{document}
\monthyear{February 2016}
\volnumber{Volume 54, Number 1}
\setcounter{page}{1}

\title{Fibonacci numbers which are products of two Pell numbers}
\author{Mahadi Ddamulira}
\address{AIMS Ghana (Biriwa)\\ 
		P.O. Box DL 676\\ 
		Adisadel, Cape Coast\\ 
		Central Region, Ghana}
\email{mahadi@aims.edu.gh}
\author{Florian Luca}
\address{School of Mathematics\\ 
		University of the Witwatersrand \\ 
		Private Bag X3, Wits 2050\\ 
		Johannesburg, South Africa}
\email{florian.luca@wits.ac.za}
\author{Mihaja Rakotomalala}
\address{AIMS Ghana (Biriwa)\\ 
		P.O. Box DL 676\\ 
		Adisadel, Cape Coast\\ 
		Central Region, Ghana}
\email{mihaja@aims.edu.gh}

\begin{abstract}
In this paper, we  find all Fibonacci numbers which are products of two Pell numbers and all Pell numbers which are products of two Fibonacci numbers.
\end{abstract} 

\maketitle

\section{Introduction}

Let $\{F_n\}_{n\ge 0} $ and $\{P_n\}_{n\ge 0}$ be the sequences of Fibonacci and Pell numbers given by $F_0=P_0=0$, $F_1=P_1=1$ and 
$$
F_{n+2}=F_{n+1}+F_n\quad {\text{\rm and}}\quad P_{n+2}= 2P_{n+1} + P_n \quad {\text{\rm for all}}\quad n\geq 0,
$$
respectively. Their first few terms are
\begin{eqnarray*}
&& \{F_n\}_{n\ge 1}\qquad 1,1,2,3,5,8,13,21,34,55,89,144,233,\ldots\\
&& \{P_n\}_{n\ge 1}\qquad 1,2,5,12,29,70,169,408,985,\ldots.
\end{eqnarray*}

 Putting $(\alpha,\beta)=((1+{\sqrt{5}})/2,(1-{\sqrt{5}})/2)$ and $(\gamma,\delta)=(1+{\sqrt{2}},1-{\sqrt{2}})$ for the pairs of roots of the characteristic equations 
$x^2-x-1=0$ and $x^2-2x-1=0$ of the Fibonacci and Pell numbers, respectively, 
then the Binet formulas for their general terms are:
$$
F_n=\frac{\alpha^n-\beta^n}{\alpha-\beta}\quad {\text{\rm and}}\quad P_n=\frac{\gamma^n-\delta^n}{\gamma-\delta}\quad {\text{\rm for~all}}\quad n\ge 0,
$$ 
respectively.

In this note, we study the Diophantine equations
\begin{equation}
\label{eq:FPP}
F_k=P_mP_n
\end{equation}
and
\begin{equation}
\label{eq:FFP}
P_k=F_m F_n,
\end{equation}
Our results are:

\begin{theorem}
\label{thm:1}
\begin{itemize}
\item[(i)] All positive integer solutions $(k,m,n)$ of equation \eqref{eq:FPP} have $k=1,2,5,12$.
\item[(ii)] All positive integer solutions $(k,\ell,m)$ of equation \eqref{eq:FFP} have $k=1,2,3,7$.
\end{itemize}
\end{theorem}

It is known that $144=12^2$ and $169=13^2$ are the largest squares in the Fibonacci and Pell sequences, respectively, and $12$ and $13$ are Pell and Fibonacci numbers, respectively. So, the above theorem says that there are no larger Fibonacci or Pell numbers which are products of two numbers from the other sequence.

When $m=1$ in equation \eqref{eq:FPP} or $k=1$ in equation \eqref{eq:FFP}, the resulting Diophantine equation is of the form
\begin{equation}
\label{eq:gen}
U_n=V_m \quad {\text{\rm for some}}\quad m,n\geq 0,
\end{equation}
where $\{U_n\}_{n\ge 0}$ and $\{V_m\}_{m\ge 0}$ are the Fibonacci and Pell sequences, respectively.  
More generally, there is a lot of literature on how to solve  equations like \eqref{eq:gen} in case $\{U_n\}_{n\ge 0}$ and $\{V_m\}_{m\ge 0}$ are  two non degenerate linearly recurrent sequences with dominant roots.  See, for example,  \cite{MM1} and \cite{MM2}. The theory of linear forms in logarithms \a la Baker gives that, under reasonable conditions (say, the dominant roots of $\{U_n\}_{n\ge 0}$ and $\{V_m\}_{m\ge 0}$ are multiplicatively independent), equation \eqref{eq:gen} has only finitely many solutions which are effectively computable. In fact, a straightforward linear form in logarithms gives some very large bounds on $\max\{m,n\}$, which then are reduced in practice either by using the LLL algorithm or by using a procedure originally discovered by Baker and Davenport \cite{BD} and perfected by Dujella and Peth\H o \cite{DP}.

\medskip

In this paper, we also use linear forms in logarithms and the Dujella-Peth\H o reduction procedure to solve equations \eqref{eq:FPP} and \eqref{eq:FFP}. 

\section{A lower bound for a linear forms in logarithms of algebraic numbers}
\label{sec:2}

In this section, we state a result concerning lower bounds for linear
forms in logarithms of algebraic numbers, which will be used in the
proof of our theorem.

Let $\eta$ be an algebraic number of degree $d$, whose
minimal polynomial over the integers is $$g(x) = a_0 \prod_{i=1}^d (x - \eta^{(i)}).$$
The logarithmic height of $\eta$ is defined as
$$
h(\eta) = \frac{1}{d}\left( \log |a_0| + \sum_{i=1}^d \log \max \{ |\eta^{(i)}|,
1\}\right).
$$
Let $\L$ be an algebraic number field and $d_{\L}$ be the degree of
the field $\L$. Let $\eta_1, \eta_2, \ldots, \eta_l \in \L$  not $0$
or $1$ and $d_1, \ldots, d_l$ be nonzero integers. We put
$$
D =\max\{|d_1|, \ldots, |d_l|, 3\},
$$
and put
$$
\Lambda = \prod_{i=1}^l \eta_i^{d_i} -1.
$$
Let $A_1, \ldots, A_l$ be positive integers such that $$A_j \geq h'(\eta_j) := \max \{d_{\L}h(\eta_j), |\log \eta_j|, 0.16
\}\quad {\text{\rm for}}\quad j=1,\ldots l.
$$
The following result is due to Matveev \cite{Matveev}.
\begin{theorem}
\label{thm:Matveev} If $\Lambda \neq 0$ and $\L \subset \R $, then
\begin{equation*}
\label{ineq:matveev} \log |\Lambda| > -1.4 \cdot
30^{l+3}l^{4.5}d_{\L}^2(1+ \log d_{\L})(1+ \log D)A_1A_2\cdots A_l.
\end{equation*}
\end{theorem}

\section{\small Proof of Theorem \ref{thm:1}}

We ran a computation for $k\le 400$ and got only the indicated solutions. We now assume that $k>400$ and that $n>m$. We do not consider the case $n=m$ since they lead to $F_k=\square$ and $P_k=\square$ whose largest solutions are $k=12$ and $k=7$, respectively, as we already pointed out in the Introduction. We deal with equation \eqref{eq:FPP} first.
We use the known inequalities that
$$
\alpha^{n-2}\le F_n\le \alpha^{n-1}\quad {\text{\rm and}}\quad \gamma^{n-2}\le P_n\le \gamma^{n-1}\quad {\text{\rm for~all}}\quad n\ge 0.
$$
Thus, 
\begin{equation}
\label{eq:bounds}
\alpha^{k-2}\le F_k=P_mP_n\le \gamma^{m+n-2}\quad {\text{\rm and}}\quad \alpha^{k-1}\ge F_k=P_nP_m\ge \gamma^{m+n-4}.
\end{equation}
Hence,
\begin{equation}
\label{eq:c1}
1+c_1(m+n-4)\le k\le 2+c_1(m+n-2),\quad {\text{\rm where}}\quad c_1=\log \gamma/\log \alpha=1.83157\ldots.
\end{equation}
In particular, $k<4n$. We get
$$
\frac{1}{\sqrt{5}}(\alpha^k-\beta^k)=\frac{1}{8}(\gamma^m-\delta^{m})(\gamma^n-\delta^n),
$$
which can be regrouped as 
$$
\left|\frac{\alpha^k}{\sqrt{5}}-\frac{\gamma^{m+n}}{8}\right|=\left|\frac{\beta^k}{\sqrt{5}}-\frac{\gamma^n \delta^m+\gamma^m\delta^n-\delta^{m+n}}{8}\right|.
$$
Since $\delta=-\gamma^{-1}$ and $\beta=-\alpha^{-1}$,  and the fact that $3/8<1/{\sqrt{5}}$, we get that
\begin{equation}
\label{eq:ttt}
\left|\frac{\alpha^k}{\sqrt{5}}-\frac{\gamma^{m+n}}{8}\right|<\frac{2}{\sqrt{5}} \max\left\{|\beta|^k, \gamma^{n-m}\right\}=\frac{2\gamma^{n-m}}{\sqrt{5}}.
\end{equation}
Dividing across by $\gamma^{m+n}/8$, we get
\begin{equation}
\label{eq:lin1}
\left|\frac{8}{\sqrt{5}} \alpha^k \gamma^{-n-m}-1\right|<\frac{16}{{\sqrt{5}} \gamma^{2m}}.
\end{equation}
On the left--hand side of \eqref{eq:lin1}  we apply Theorem \ref{thm:Matveev} with the data 
$$
l=3,~\eta_1=8/{\sqrt{5}},~\eta_2=\alpha,~\eta_3=\gamma,~d_1=1,~d_2=k,~d_3=-m-n.
$$ 
We take ${\mathbb L}={\mathbb Q}({\sqrt{2}}, {\sqrt{5}})$, for which $d_{\mathbb L}=4$. Since 
$$
h(\eta_1)=\log 8,~h(\eta_2)=(1/2)\log \alpha,~h(\eta_3)=(1/2)\log \gamma,
$$ 
we take $A_1=4\log 8,~A_2=2\log \alpha,~A_3=2\log \gamma$. Finally, we can take $D=4n$. Note that
$$
\Lambda_1=\frac{8}{\sqrt{5}} \alpha^k \gamma^{-n-m}-1.
$$
The fact that it isn't zero follows from the fact that if it were, we would then get that $\alpha^{-k}\gamma^{m+n}=8/{\sqrt{5}}$. However, the left-hand side of the above relation is a unit in ${\mathbb L}$, whereas the right hand side is not as its norm over ${\mathbb L}$ is $2^{12}/5^2$. Thus, $\Lambda_1\ne 0$. Theorem \ref{thm:Matveev} gives that 
$$
\log |\Lambda_1|>-1.4\times 30^{6} \times 3^{4.5} \times 4^2 (1+\log 4)(1+\log (4n)) (4\log 8 ) (2\log(\alpha)) (2\log \gamma).
$$
Comparing the above inequality with \eqref{eq:lin1}, we get
\begin{equation}
\label{eq:2m}
2m\log \gamma-\log(16/{\sqrt{5}})<7.8\times 10^{13} (1+\log (4n)).
\end{equation}
Hence,
\begin{equation}
\label{eq:boundfor2m}
m\log \gamma<4\times 10^{13} (1+\log(4n)).
\end{equation}
Next we return to equation \eqref{eq:FPP} and rewrite it as 
$$
\left|\frac{\alpha^k}{{\sqrt{5}} P_m} -\frac{\gamma^n}{2{\sqrt{2}}}\right|=\left|\frac{\beta^k}{{\sqrt{5}}P_m}-\frac{\delta^n}{2{\sqrt{2}}}\right|\le \frac{2}{\sqrt{5}} \max\left\{\frac{1}{\alpha^k}, \frac{1}{\gamma^n}\right\}.
$$
We divide both sides above by $\gamma^n/2{\sqrt{2}}$ getting
$$
\left|\frac{2{\sqrt{2}}}{{\sqrt{5}}P_m} \alpha^k \gamma^{-n}-1\right|\le \frac{4{\sqrt{2}}}{\sqrt{5}} \max\left\{\frac{1}{\alpha^k \gamma^n}, \frac{1}{\gamma^{2n}}\right\}.
$$ 
From \eqref{eq:bounds}, we get that 
$$
\frac{1}{\alpha^k\gamma^n}= \frac{1/\alpha}{\alpha^{k-1}\gamma^n}\le \frac{(1/\alpha)}{\gamma^{2n+m-4}}
=\frac{\gamma^3/\alpha}{\gamma^{2n+m-1}}<\frac{9}{\gamma^{2n}},
$$
because $\gamma^3/\alpha<9$ and $m\ge 1$. Thus,
\begin{equation}
\label{eq:lin2}
\left|\frac{2{\sqrt{2}}}{{\sqrt{5}}P_m} \alpha^k \gamma^{-n}-1\right|\le \frac{4{\sqrt{2}}\times 9}{{\sqrt{5}}\gamma^{2n}}=\frac{36{\sqrt{2}}}{{\sqrt{5}}\gamma^{2n}}.
\end{equation}
On the left--hand side of \eqref{eq:lin2}  we apply Theorem \ref{thm:Matveev} with the data 
$$
l=3,~\eta_1={\sqrt{5}}P_m/2{\sqrt{2}},~\eta_2=\alpha,~\eta_3=\gamma,~d_1=-1,~d_2=k,~d_3=-n.
$$ 
We take again ${\mathbb L}={\mathbb Q}({\sqrt{2}}, {\sqrt{5}})$, for which $d_{\mathbb L}=4$. As before,
$$
h(\eta_2)=(1/2)\log \alpha,~h(\eta_3)=(1/2)\log \gamma,
$$ 
so we can take $A_2=2\log \alpha,~A_3=2\log \gamma$. As for $h(\eta_1)$, the polynomial 
$$
8X^2-5 P_m^2
$$
has $\eta_1$ as a root. Thus, 
\begin{eqnarray*}
h(\eta_1) & \le & \frac{1}{2}\left (\log 8+2\log ({\sqrt{5}} P_m/2{\sqrt{2}})\right)\\
& = & \log P_m+\log{\sqrt{5}}\le (m-1)\log \gamma+\log {\sqrt{5}}\\
& < & m\log \gamma.
\end{eqnarray*}
Using \eqref{eq:boundfor2m}, we can take
$$
A_1= 16\times 10^{13} (1+\log(4n))>4 h(\eta_1).
$$
Finally, we can take $B=4n$. Note that
$$
\Lambda_2=\frac{2{\sqrt{2}}}{{\sqrt{5}} P_m} \alpha^k \gamma^{-n-m}-1.
$$
Similarly to the argument used to prove that $\Lambda_1\ne 0$, one justifies that $\Lambda_2\ne 0$. Theorem \ref{thm:Matveev} gives that 
$$
\log |\Lambda_2|>-1.4\times 30^{6} \times 3^{4.5} \times 4^2 (1+\log 4)(1+\log (4n))^2 16 \times 10^{13}\times  (2\log(\alpha)) (2\log \gamma).
$$
Comparing this with \eqref{eq:lin2}, we get
$$
2n\log \gamma-\log(36 {\sqrt{2}}/{\sqrt{5}})<1.5\times 10^{27} (1+\log 4n))^2,
$$
giving 
\begin{equation}
\label{eq:5}
n<5\times 10^{30}. 
\end{equation}
The same arguments apply to equation \eqref{eq:FFP} (just swap the roles of the pairs $(\alpha,\beta)$ and $(\gamma,\delta)$ of $1/{\sqrt{5}}$ and $1/(2{\sqrt{2}})$. Let us give the details. We assume $m\ge 3$, otherwise $m\in \{1,2\}$, $F_m=1$ and the solutions of \eqref{eq:FFP} are among the solutions to \eqref{eq:FPP} with $m=1$. Inequality 
\eqref{eq:c1} becomes
\begin{equation}
\label{eq:c2}
1+c_2(m+n-4)\le k\le 2+c_2(m+n-2),\qquad c_2=1/c_1=\log \alpha/\log \gamma=0.545979\ldots,
\end{equation}
which implies in particular that $k\le 3n$. The analog of inequality \eqref{eq:ttt} is 
\begin{eqnarray}
\label{eq:ttt1}
\left|\frac{\gamma^k}{2\sqrt{2}}-\frac{\alpha^{m+n}}{5}\right| & = & \left|\frac{\delta^k}{2{\sqrt{2}}}-\frac{\alpha^{n}\beta^m+\alpha^m \beta^n-\beta^{m+n}}{5}\right|\\
& \le & \frac{6}{5} \max\left\{|\delta|^k, \alpha^{n-m}\right\}=\frac{6\alpha^{n-m}}{5}.
\end{eqnarray}
This leads to
\begin{equation}
\label{eq:lin3}
\left|\frac{5}{2{\sqrt{2}}}\gamma^k \alpha^{-n-m}-1\right|<\frac{6}{\alpha^{2m}},
\end{equation}
which is the analogue of \eqref{eq:lin1}. We check that the amount  $\Lambda_3$ in the left--hand side above is non-zero by an argument similar to the one used to prove that $\Lambda_1$ and $\Lambda_2$ are non-zero, and apply Theorem \ref{thm:Matveev} to get a lower bound for it,
getting 
$$
\log \Lambda_3>-1.4\times 30^{6} \times 3^{4.5} \times 4^2 (1+\log 4) (1+\log (3n)) (4\log 5) (2\log(\alpha)) (2\log \gamma).
$$
We get that the analog of \eqref{eq:2m} is
$$
2m\log \alpha-\log 6<5.98\times 10^{13}(1+\log(3n)),
$$
giving 
\begin{equation}
\label{eq:uuu}
m\log\alpha+1<3\times 10^{13} (1+\log(3n)),
\end{equation}
which is the analog of inequality \eqref{eq:boundfor2m}. Returning to equation \eqref{eq:FFP}, we get
\begin{equation}
\label{eq:10}
\left|\frac{\gamma^k}{2{\sqrt{2}}F_m}-\frac{\alpha^n}{\sqrt{5}}\right|=\left|\frac{\delta^k}{2{\sqrt{2}}F_m}-\frac{\beta^n}{\sqrt{5}}\right|\le \frac{2}{\sqrt{5}} \max\left\{\frac{1}{\gamma^k},\frac{1}{\alpha^n}\right\}.
\end{equation}
By \eqref{eq:c2}, we get
$$
\gamma^k\ge \gamma \alpha^{m+n-4}\ge \gamma \alpha^{-3} \alpha^n,
$$
so
\begin{equation}
\label{eq:11}
\frac{1}{\gamma^k}\le \frac{\alpha^3/\gamma}{\alpha^n}<\frac{2}{\alpha^n}.
\end{equation}
Hence, by \eqref{eq:10} and \eqref{eq:11}, we get
\begin{equation}
\label{eq:lin4}
\left|\frac{\sqrt{5}}{2{\sqrt{2}}F_m} \gamma^k \alpha^{-n}-1\right|<\frac{4}{\alpha^{2n}}.
\end{equation}
This is the analog of \eqref{eq:lin2}. Writing $\Lambda_4$ for the amount under the absolute value in the left--hand side above, we get that it is not $0$ by arguments similar to the ones used to prove that 
$\Lambda_i\ne 0$ for $i=1,2,3$. We apply Matveev's theorem as we did for $\Lambda_2$. Here, $\eta_1=2{\sqrt{2}} F_m/{\sqrt{5}}$ is a root of $5X^2-8F_m^2$. Its height therefore satisfies 
\begin{eqnarray*}
h(\eta_1) & \le & \log F_m+\log 2{\sqrt{2}}\le (m-1)\log \alpha+\log 2{\sqrt{2}}\\
& < & m\log \alpha+1<3\times 10^{13} (1+\log(3n)),
\end{eqnarray*}
by \eqref{eq:uuu}. We get that 
$$
\log \Lambda_4>-1.4\times 30^{6} \times 3^{4.5} \times 4^2 (1+\log 4)(1+\log (3n))^2 12 \times 10^{13}\times  (2\log(\alpha)) (2\log \gamma),
$$
which together with \eqref{eq:lin4} leads to 
$$
2n\log \alpha-\log 4<1.2\times 10^{27} (1+\log(3n))^2,
$$
giving
$$
n<7\times 10^{30}.
$$
So, comparing the above bound with \eqref{eq:5}, we conclude that both in equation \eqref{eq:FPP} and \eqref{eq:FFP}, we get $n<7\times 10^{30}$.  We record what we proved as a lemma.

\begin{lemma}
If $(k,m,n)$ are positive integers satisfying one of the equations \eqref{eq:FPP} or \eqref{eq:FFP} with $m\le n$, then $k<4n$ and $n<7\times 10^{30}$.  
\end{lemma}

Now we need to reduce the bound. To do so, we make use several times of the following result, which is a slight variation of a result due to Dujella and Peth\H{o} \cite{DP} which itself is a generalization of a result of Baker and Davenport \cite{BD}. 
The proof is almost identical to the proof of the corresponding result in \cite{DP} and the details have been worked out in Lemma 2.9 in \cite{THesisJJ}. For a real number $x$, we put  $ ||x|| = \min\{|x-n|\,:\, n \in \Z\}$ for the distance from $x$ to the nearest integer.
\begin{lemma}
\label{Dujella-Petho}
Let $M$ be a positive integer, let $p/q$ be a convergent of the continued
fraction of the irrational
$\tau$ such that $q > 6M$, and let $A, B,  \mu$ be some real
numbers with $A > 0$ and $B > 1$. Let $\epsilon :=||\mu q|| - M||\tau q||$. If $\epsilon > 0$, then there is no solution to the
inequality
$$
0 < m\tau - n + \mu < AB^{-k},
$$
in positive integers $m, n$ and $k$ with
$$
m \leq M \qquad {\text{and}} \qquad k \geq \dfrac{\log(Aq/\epsilon)}{\log B}.
$$
\end{lemma}

We look at \eqref{eq:lin1}. Assume that $m\ge 20$. Put 
$$
\Gamma_1:=k\log \alpha-(n+m)\log \gamma+\log(8/{\sqrt{5}}).
$$
Then $|e^{\Gamma_1}-1|=|\Lambda_1|<1/4$ by \eqref{eq:lin1}, which implies that $|\Gamma_1|<1/2$. 
Since $|x|<2|e^{x}-1|$ whenever $x\in (-1/2,1/2)$, we get from $\Lambda_1=e^{\Gamma_1}$ and \eqref{eq:lin1} that 
$$
|\Gamma_1|<\frac{32}{{\sqrt{5}}\gamma^{2m}}.
$$
If $\Gamma_1>0$, then 
$$
0<k\left(\frac{\log \alpha}{\log \gamma}\right)-(n+m)+\frac{\log(8/{\sqrt{5}})}{\log \gamma}<\frac{32}{({\sqrt{5}} \log \gamma) \gamma^{2m}}<\frac{17}{\gamma^{2m}}.
$$
We apply Lemma \ref{Dujella-Petho} with $M=3\times 10^{31}$ (note that $M>4n>k$), 
$$
\tau=\frac{\log \alpha}{\log \gamma},\quad \mu=\frac{\log(8/{\sqrt{5}})}{\log\delta},\quad A=17,\quad B=\gamma^2.
$$
Writing $\tau=[a_0,a_1,\ldots]$ as a continued fraction, we get 
$$
[a_0,\ldots,a_{74}]=\frac{p_{74}}{q_{74}}=\frac{2037068391552562960855777461929676271}{
3731035235978315437343082205475618926},
$$
and we get
$q_{74}>3\times 10^{36}>6M$. We compute $\varepsilon=\| \mu q_{74}\|-M\| \tau q_{74}\|>0.4.$ The reason that we picked the $74$th convergent is that both the inequalities $q_{74}>6M$ and $\varepsilon>0$ hold. 
 Thus, by Lemma \ref{Dujella-Petho}, we get
$m\le 49$. A similar conclusion is reached if we assume that $\Gamma_1<0$. This was in the case of inequality \eqref{eq:lin1}. In the case of inequality \eqref{eq:lin3}, assuming again 
that $m\ge 20$, we get that 
$$
\left|(n+m)\log \alpha-k\log \gamma-\log(5/2{\sqrt{2}})\right|<\frac{12}{5\alpha^{2m}}.
$$
Let $\Gamma_3$ be the expression under the absolute value of the left--hand side above. If $\Gamma_3>0$, we get
$$
0<(n+m)\left(\frac{\log \alpha}{\log \gamma}\right)-k+\frac{\log(2{\sqrt{2}}/5)}{\log \gamma}<\frac{12}{(5\log \gamma) \alpha^{2m}}<\frac{3}{\alpha^{2m}}.
$$
We keep the same values for $M,~\tau,~q$ and only change $\mu$ to 
$$
\mu'=\frac{\log(2{\sqrt{2}}/5)}{\log \gamma},\quad A=3,\quad B=\alpha^2.
$$
We get $\varepsilon>0.2$, and by Lemma \ref{Dujella-Petho}, $m\le 90$.
A similar conclusion is reached if $\Gamma_3<0$.  Thus, $m\le 90$ in all cases. Now we move on to \eqref{eq:lin2}. Assume $n>100$. We then get 
$$ 
\left|k\log \alpha-n\log \gamma+\log(2{\sqrt{2}}/{\sqrt{5}}P_m)\right|<\frac{72{\sqrt{2}}}{{\sqrt{5}}\delta^{2n}}.
$$
Let $\Gamma_2$ be the expression under the absolute value in the left--hand side above. If $\Gamma_2>0$, we then get
$$
0<k\left(\frac{\log \alpha}{\log \gamma}\right)-n+\frac{\log(2{\sqrt{2}}/({\sqrt{5}}P_m))}{\log \gamma}<\frac{72 {\sqrt{2}}}{({\sqrt{5}} \log \gamma)\gamma^{2n}}<\frac{52}{\gamma^{2n}}.
$$
We keep the same values for $M,~\tau,~q$ and only change $\mu$ to 
$$
\mu_m=\frac{\log(2{\sqrt{2}}/({\sqrt{5}} P_m))}{\log \gamma},\quad A=52,\quad B=\gamma^2\quad {\text{\rm for~all}}\quad m=1,\ldots,90.
$$
We get $\varepsilon>0.019$, so $n\le 53$. A similar conclusion is reached if $\Gamma_2<0$.  Finally, if instead of \eqref{eq:lin2}, we have \eqref{eq:lin4}, then a similar argument leads to 
$$
\left|n\log \alpha- k\log \gamma+\log(2{\sqrt{2}}F_m/{\sqrt{5}})\right|<\frac{4}{\alpha^{2n}}.
$$
Putting $\Gamma_4$ for the amount under the absolute value in the left--hand side above, we get in case $\Gamma_4>0$ that
$$
0<n \left(\frac{\log \alpha}{\log \gamma}\right)-k+\frac{\log(2{\sqrt{2}} F_m/{\sqrt{5}})}{\log \gamma}<\frac{4}{\log \gamma \alpha^{2n}}<\frac{5}{\alpha^{2n}}.
$$
We keep the same values for $M,~\tau,~q$ and only change $\mu$ to 
$$
\mu_m=\frac{\log(2{\sqrt{2}} F_m/{\sqrt{5}})}{\log \gamma},\quad A=5,\quad B=\alpha^2,\quad {\text{\rm for~all}}\quad m=1,\ldots,90.
$$
We get $\varepsilon>0.005$, so $n\le 94$. So, in all cases $n\le 94$, so $k<400$.  We generated $\{F_k\}_{1\le k\le 400}$ and $\{P_m P_n\}_{1\le m<n\le 100}$ and intersected them,   and also
$\{P_k\}_{1\le k\le 400}$ and $\{F_mF_n\}_{1\le m<n\le 100}$ and intersected them and got no other solutions. Hence, Theorem \ref{thm:1} is proved.

\section{Comments}

It is apparent from our proof that the method is more general and shows that every equation of the form
$$
U_k=V_m V_n
$$
has only finitely effectively computable many positive integer solutions $(k,m,n)$ provided that $\{U_n\}_{n\ge 0}$ and $\{V_n\}_{n\ge 0}$ satisfy a few technical conditions such as:
\begin{itemize}
\item[(i)] they are both non degenerate binary recurrent and have characteristic equations of real roots $\alpha,~\beta$ and $\delta,~\gamma$ with $\alpha\beta=\pm 1$ and $\gamma\delta=\pm 1$.
\item[(ii)] ${\mathbb Q}[\alpha]$ and ${\mathbb Q}[\delta]$ are distinct quadratic fields.
\end{itemize}
In fact, more is true, namely that for fixed $k$ and  $s$, the diophantine equation
$$
\prod_{i=1}^k F_{n_i}=\prod_{j=1}^s P_{m_j}
$$
has only finitely many positive integer solutions 
$$
(n_1,\ldots,n_k,m_1,\ldots,m_k)
$$ 
and all such are effectively computable. 
Such a statement is not very difficult to prove. A deeper conjecture made in \cite{LPW} to the effect that the intersection of the multiplicative group generated by $\{F_n\}_{n\ge 1}$ with the multiplicative group generated by Pell numbers $\{P_n\}_{n\ge 1}$ is finitely generated cannot unfortunately be attacked by these methods. 

\section*{Acknowledgements} We thank the referee for comments which improved the quality of this manuscript.

\medskip

\noindent MSC2010: 11B39, 11D61.


\begin{thebibliography}{99}

\bibitem{BD} A. Baker and H. Davenport, \emph{The equations $3x^2-2=y^2$ and $8x^2-7=z^2$}, Quart.J of Math. Ser.(2) \textbf{20} (1969), 129--137.

\bibitem{THesisJJ} J. J. Bravo, \emph{Arithmetic properties of generalized Fibonacci sequences}, Ph. D. Thesis, Universidad Nacional Aut\'onoma de M\'exico, 2014.  

\bibitem{DP} A.~Dujella and A.~Peth\H o, \emph{A generalization of a theorem of
Baker and Davenport}, Quart. J. Math. Oxford Ser. (2) \textbf{49} (1998), 291--306.

\bibitem{LPW} F. Luca, C. Pomerance and S. Wagner, \emph{On Fibonacci Integers}, J. Number Theory \textbf{131} (2011), 440--457.

\bibitem{Matveev} E. M. Matveev, \emph{An explicit lower bound for a homogeneous rational linear form in the logarithms of algebraic numbers II}, Izv. Ross. Akad. Nauk Ser. Mat.  \textbf{64} (2000), 125--180; translation in Izv. Math. \textbf{64} (2000), 1217--1269.


\bibitem{MM1} M. Mignotte, \emph{Intersection des images de certaines suites r\'eccurentes lin\'eaires}, Theor. Comp. Sci. \textbf{7} (1978), 117--121.

\bibitem{MM2} M. Mignotte, \emph{Une extention du th\'eoreme de Skolem-Mahler},  C. R. Acad. Sci. Paris \textbf{288} (1979), 233--235.


\end{thebibliography}
\end{document}